\input amssym.def
\input amssym.tex
\newtheorem{cast}{Claim.}
\newtheorem{claim:q+1}[cast]{Claim.}
\newenvironment{proof}{\begin{trivlist}\item[\hskip%
\labelsep{\bf Proof.\quad}]}%
{\hfill{$\Box$}\end{trivlist}}
\newcommand{\casa}{Author:
{\tt$\backslash$def$\backslash$casa$\{$??$\}$}}
\newcommand{\casatwo}{}
\newcommand{\finito}{
\par\vspace*{1cm}\noindent
\begin{minipage}[t]{.4\hsize}
\baselineskip=11pt\parindent0pt\noindent\casa\hfil
\end{minipage}\hfill
\begin{minipage}[t]{.4\hsize}
\baselineskip=11pt\parindent0pt\noindent\casatwo\hfil
\end{minipage}
\hfill\break
\vskip.4truecm}

\def\d{{\cal D}}

\def\div{{\rm div}_\infty}
\documentstyle[12pt]{article}
\renewcommand{\casa}{Fachbereich 6

Mathematik und Informatik

Universit\"at Essen

W - 4300 Essen 1, FRG

}
\renewcommand{\casatwo}{Mathematics Section

ICTP

P.O. Box 586 - 34100

Trieste - Italy

}
\begin{document}
\title{A note on the genus of certain curves\\
 over finite fields}
\author{Rainer Furhmann \and Fernando Torres\thanks{Supported by a grant
from the International Atomic Energy Agency and UNESCO.}}
\date{}
\maketitle
\begin{abstract}
We prove the following result which was conjectured by Stichtenoth and
Xing: let $g$ be the genus of a projective, irreducible non-singular curve
over the finite field $\Bbb F_{q^2}$ and whose number of
$\Bbb F_{q^2}$-rational points attains the Hasse-Weil bound; then either $4g\le
(q-1)^2$ or $2g=(q-1)q$.
\end{abstract}
\noindent {\bf 1.} Throughout this note by a curve we mean a projective,
irreducible non-singular algebraic curve defined over the finite field
$k:= \Bbb F_{q^2}$ with $q^2$ elements. This note is concerning with the
genus $g$ of a curve $X$ whose number
of $k$-rational points $\#X(k)$ attains the Hasse-Weil bound, that is,
\begin{equation}\label{max}
\#X(k)=q^2+2gq+1.
\end{equation}
These curves are very useful for applications to coding theory
(see e.g. \cite{Sti}, \cite{Tsf-Vla}). It is known that $2g\le (q-1)q$
(\cite[V.3.3]{Sti}), and  that the Hermitian curve is the only one
that satisfies (\ref{max}) and whose genus is $(q-1)q/2$
\cite{R-Sti}. Here we prove the following result which was conjectured by
Stichtenoth and Xing in \cite{Sti-X}.
\medskip

\noindent {\bf Theorem.}\quad Let $X$ be a curve of genus $g$, and
suppose that $X$ satisfies (\ref{max}). Then
$$
4g\le (q-1)^2\qquad {\rm or}\qquad 2g=(q-1)q.
$$

We prove this theorem by using \cite[Prop. 1]{Sti-X}, R\"uck and
Stichtenoth's \cite[Lemma 1]{R-Sti} and by using a particular case of the
the theory of Frobenius orders sequence associated to linear system of
curves which was developed by St\"ohr and Voloch \cite{S-V} as a way of
improving the Hasse-Weil bound. We also state another proof of
R\"uck and Stichtenoth's result (cf. \cite{R-Sti}) concerning Hermitian
curves (see \S4). We remark that Hirschfeld, Storme, Thas and Voloch also
stated a characterization of Hermitian curves by using some
results from \cite{S-V} (see
\cite{HSTV}). In \S3 we use ideas from the proof of \cite[Lemma 1]{HSTV}.

In \cite{F-T} is considered the case of curves satisfying (\ref{max}) and
whose genus is bounded from above by $(q-1)^2/4$.

We are indebted to Professor J.F. Voloch for pointing out to us that the
proof of the theorem in the previous version of this paper was incomplete.
\bigskip

\noindent {\bf 2.} Let $X$ be a curve of genus $g$, satisfying the
hypothesis (\ref{max}). The starting point of the proof is the fact that
there exists a $k$-rational point $P_0\in X$ such that $q$ and $q+1$ are
not-gaps at $P_0$ \cite[Prop. 1]{Sti-X}. Then the linear system
$$
\d=g^{r+1}_{q+1}:=|(q+1)P_0|
$$
is simple. Thus the genus $g$ of $X$ and the dimension $r+1$ of $\d$ are
related by
Castelnuovo's genus bound for curves in projective spaces (\cite{C},
\cite[p. 116]{ACGH}, \cite[Corollary 2.8]{Ra}). Thus we have
\begin{equation}\label{cast}
2g\le M(q-r+e),
\end{equation}
where $M$ is the biggest integer $\le q/r$ and $e:= q-Mr$.
\begin{cast}[[Sti-X, Prop. 3{]}]
If $r+1\ge 3$, then $4g\le (q-1)^2$.
\end{cast}
\begin{proof}
{}From (\ref{cast}) we have
$$
2g \le {(q-e)(q-r+e)\over r}\le {(2q-r)^2 \over4r}.
$$
This implies the result.
\end{proof}

{}From now on we considerer $r+1=2$. We apply \cite{S-V} to the
linear system $\d = g^2_{q+1} = |(q+1)P_0|$.
For $P\in X$ let $0=j_0(P)<j_1(P)<j_2(P)\le q+1$ denote the $(\d,
P)$-orders (\cite[\S1]{S-V}. The fact that $q$ and $q+1$ are non-gaps at
$P_0$ implies
$$
j_1(P_0)=1\qquad {\rm and}\qquad j_2(P_0)=q+1.
$$
To compute $j_2(P)$ for $P\in X(k)$ we use the fact that the linear
system $|(q+1)P|$ is equivalent to $\d$ (\cite[Lemma 1]{R-Sti}). Then
\begin{equation}\label{j2}
j_2(P)=q+1,
\end{equation}
provided $P\in X(k)$.

Let denote by $0=\nu_0<\nu_1$ the $k$-Frobenius order sequence of $\d$
(see \cite[Prop. 2.1]{S-V}). $\nu_1$ satisfies
\begin{equation}\label{desnu}
\nu_1\le j_2(P)-j_1(P)=q+1-j_1(P),
\end{equation}
for each $P\in X(k)$ \cite[Corollary 2.6]{S-V}.

Associated to $\d$ there exists a positive
divisor $S$ of degree
\begin{equation}\label{degS}
{\rm deg}(S)=\nu_1(2g-2)+(q^2+2)(q+1)
\end{equation}
whose support contains $X(k)$. Moreover by \cite[Prop. 2.4]{S-V} we have
\begin{equation}\label{maS}
v_P(S)\ge j_1(P)+(j_2(P)-\nu_1)\ge q+2-\nu_1.
\end{equation}
\begin{claim:q+1}\label{:q+1} Suppose $4g>(q-1)^2$. Then
\begin{list}
\setlenght{\rightmargin 0cm}{\leftmargin 0cm}
\itemsep=0.5pt
\item[(i)] \quad $\nu_1=q$.
\item[(ii)] \quad $j_1(Q)=1$ for $Q\in X(k)$.
\end{list}
\end{claim:q+1}
\begin{proof}
Statement (ii) follows from statement (i) and relation (\ref{desnu}).
Applying  (\ref{desnu}) for $P_0$ we also obtain $\nu_1\le q$.
Now by using relations (\ref{maS}), (\ref{max}) and (\ref{degS})
we get
$$
(q-1)(\nu_1(q+1)-q)\ge 2g(q^2-\nu_1(q+1)+2q).
$$
Finally by using the hypothesis on $g$ we get $\nu_1\ge q$ and
we are done.
\end{proof}
\bigskip

\noindent {\bf 3. Proof of the Theorem.}\quad Let assume that
$$
4g> (q-1)^2.
$$
{}From the fact that $q$ and $q+1$ are non-gaps at $P_0$ we have $2g\le
(q-1)q$ (see e.g Jenkins's \cite{J}).

Let $R$ be the
divisor on $X$ supporting the $\d$-Weierstrass points,
and let denote by $0=\epsilon_0<\epsilon_1<\epsilon_2$ the order
sequence of $\d$ (see \cite[\S1]{S-V}). The degree of $R$ is given by
$$
{\rm deg}(R)= (\epsilon_1+\epsilon_2)(2g-2)+3(q+1).
$$
For each $P\in X$ we have $v_P(R)\ge \sum_{i=0}^{2}(j_i(P)-\epsilon_i)$.
Now since $\epsilon_i\le j_i(P)$ for
each $i$, from Claim \ref{:q+1} we have $\epsilon_1=1$. Moreover by
\cite[Prop. 2.1]{S-V} we have $\epsilon_2=\nu_1$. Thus $\epsilon_2=q$ by
Claim \ref{:q+1}. Hence by relation (\ref{j2}) we have
$$
v_P(R)\ge 1
$$
for each $P\in X(k)$. Consequently by $(\ref{max})$ and ${\rm deg}(R)$ it
follows that $2g\ge (q-1)q$, which proves the result.\hfill{$\Box$}
\bigskip

\noindent {\bf 4. Remark.}\quad We close this note by proving
that a curve satisfying (\ref{max}) and whose genus is $(q-1)q/2$ is
$k$-isomorphic to the so called Hermitian curve which is defined by an
equation of type
$$
y^q+y=x^{q+1}.
$$
This result
was proved in \cite{R-Sti}. Here we give a different proof which
is inspired on the example stated in \cite[p. 16]{S-V}.

Let $P_0$ be as in \S2 and let $x, y \in k(X)$ such that
$$
\div(x)=qP_0\qquad {\rm and}\qquad \div(y)=(q+1)P_0.
$$
Then by the Riemann-Roch theorem the $k$-dimension of the
linear system $|q(q+1)P_0|$ is equal to $(q+1)(q+2)/2 -1$. Since
$$
\#\{(i,j) \in \Bbb N \times \Bbb N: iq+j(q+1) \le q(q+1) \} =
{(q+1)(q+2)\over 2}+1,
$$
there exists a non-trivial $k$-linear relation:
\begin{equation}\label{rel0}
F=\mathop{\sum}\limits^{}_{iq+j(q+1)\le q(q+1)}{a_{i,j}x^iy^j} = 0,
\end{equation}
where $a_{q+1,0} \not= 0$ and $a_{0,q}\not = 0$. $X$ is isomorphic to the
plane curve defined by $F$. After a $k$-projective transformation we may
assume $a_{q+1,0}=a_{0,q}=1$ and $a_{q,0}=0$. The fact that $\nu_1=q>1$
implies
\begin{equation}\label{rel}
(x-x^{q^2})F_x=(y^{q^2}-y)F_y,
\end{equation}
where $F_x$ (resp. $F_y$) stands for the partial derivative
with respect to the variable $x$ (resp. $y$) of $F$. Let $v$ denote
the valuation associated to $P_0$. From (\ref{rel}) we find that
$F_y\not=0$ and $v(F_y)=0$. Thus $F_y=a_{0,1}$. Now we compute $y^{q^2}$
be means of (\ref{rel0}) and hence from (\ref{rel}) we get $a_{0,1}=1$,
$a_{i,j}=0$ for $(i,j)\not\in\{(0,0)), (0,q), (q+1,0)\}$. Finally by
another $k$-projective transformation we may assume $a_{0,0}=0$ and we
find that $F$ defines the Hermitian curve.

\finito

\begin{thebibliography}{Feto 99}

\bibitem[ACGH]{ACGH} Arbarello, E.; Cornalba, M.; Griffiths, P.A. and
Harris, J.: Geometry of algebraic curves, Vol. I, Springer-Verlag, New
York 1985.

\bibitem[C]{C} Castelnuovo, G.: Ricerche di geometria sulle curve
algebriche, Atti. R. Acad. Sci. Torino {\bf 24}, 196--223 (1889).

\bibitem[F-T]{F-T} Furhmann, R.; Torres, F.: Curves over finite fields
with maximal number of rational points, (In preparation).

\bibitem[HSTV]{HSTV} Hirschfeld, J.W.P.; Storme, L.; Thas, J.A. and
Voloch, J.F.: A characterization of Hermitian curves, J. Geom. {\bf 41},
72--78 (1991).

\bibitem[J]{J} Jenkins, J.A.: Some remarks on Weierstrass points, Proc.
Amer. Math. Soc. {\bf 44}, 121--122 (1974).

\bibitem[Ra]{Ra} Rathmann, J.: The uniform position principle for curves
in characteristic $p$, Math. Ann. {\bf 276}, 565--579 (1987).

\bibitem[R-Sti]{R-Sti} R\"uck, H.G.; Stichtenoth, H.: A characterization of
Hermitian function
fields over finite fields, J. Reine. Angew. Math. {\bf 457}, 185--188 (1994).

\bibitem[Sti]{Sti} Stichtenoth, H.: Algebraic functions fields and codes,
Springer-Verlag, Berlin, 1993.

\bibitem[Sti-X]{Sti-X} Stichtenoth, H.; Xing, C.: The genus of maximal
functions fields, Manuscripta Math. {\bf 86}, 217--224 (1995).

\bibitem[S-V]{S-V} St\"ohr, K.O.; Voloch, J.F.: Weierstrass points and
curves over finite fields, Proc. London Math. Soc. (3) {\bf 52}, 1--19 (1986).

\bibitem[Tsf-Vla]{Tsf-Vla} Tsfasman, M.; Vladut, S.G.:
Algebraic-Geometric Codes, Kluwer Academic Publishers,
Dordrecht-Boston-London 1991.

\end{thebibliography}
\end{document}